\documentclass{article}

\usepackage{PRIMEarxiv}

\usepackage[utf8]{inputenc} 
\usepackage[T1]{fontenc}    
\usepackage{hyperref}       
\usepackage{url}            
\usepackage{booktabs}       
\usepackage{amsfonts}       
\usepackage{nicefrac}       
\usepackage{microtype}      
\usepackage{lipsum}
\usepackage{fancyhdr}       
\usepackage{graphicx}       
\graphicspath{{media/}}     
\usepackage{xcolor}
\usepackage{colortbl}
\usepackage{subcaption}
\newtheorem{definition}{Definition}%
\usepackage{wrapfig}
\usepackage{multirow}
\usepackage{graphicx}
\usepackage{float}
\title{Hypergraph Classification via Persistent Homology}

\author{
  Mehmet Emin Aktas \\
  Institute for Insight \\
  Georgia State University \\
  Atlanta, GA\\
  \texttt{maktas@gsu.edu} \\
   \And
   Thu Nguyen, Rakin Riza \\
  Department of Computer Science \\
  University of Central Oklahoma \\
  Edmond, OK\\
  \texttt{\{tnguyenthuyanh, rriza\}uco@edu} \\
  \And
   Muhammad Ifte Islam, Esra Akbas \\
  Department of Computer Science \\
 Georgia State University \\
  Atlanta, GA\\
  \texttt{\{mislam29, eakbas1\}gsu@edu} \\
}

\begin{document}
\maketitle

\begin{abstract}
Persistent homology is a mathematical tool used for studying the shape of data by extracting its topological features. It has gained popularity in network science due to its applicability in various network mining problems, including clustering, graph classification, and graph neural networks. The definition of persistent homology for graphs is relatively straightforward, as graphs possess distinct intrinsic distances and a simplicial complex structure. However, hypergraphs present a challenge in preserving topological information since they may not have a simplicial complex structure. In this paper, we define several topological characterizations of hypergraphs in defining hypergraph persistent homology to prioritize different higher-order structures within hypergraphs. We further use these persistent homology filtrations in classifying four different real-world hypergraphs and compare their performance to the state-of-the-art graph neural network models. Experimental results demonstrate that persistent homology filtrations are effective in classifying hypergraphs and outperform the baseline models. To the best of our knowledge, this study represents the first systematic attempt to tackle the hypergraph classification problem using persistent homology.
\end{abstract}


\section{Introduction}

Hypergraphs have emerged as a powerful tool for modeling complex systems that involve higher-order interactions, where entities are depicted as nodes, and higher-order interactions are represented as hyperedges~\cite{singh2022survey,bretto2013applications}. This approach provides a more accurate representation of complex systems than graphs, as graphs are limited to capturing pairwise relationships between entities. In contrast, complex systems require the description of relations as higher-order interactions rather than just pairwise relationships~\cite{battiston2020networks}. For example, in coauthorship networks, authors are represented as nodes, and articles are represented as hyperedges. Similarly, in drug-drug interaction networks, substances that constitute drugs are represented as nodes, and drugs themselves are represented as hyperedges. Hypergraphs have numerous applications in network science, including but not limited to music recommendation~\cite{bu2010music}, social relationship analysis~\cite{yang2019revisiting,feng2019hypergraph}, gene expression classification~\cite{tian2009hypergraph}, and drug-drug interaction detection~\cite{saifuddin2022hygnn}.

Furthermore, persistent homology (PH) has gained widespread popularity as a powerful tool for extracting the topological characteristics of graphs and solving various graph mining problems~\cite{aktas2019persistence}. The fundamental concept of PH involves transforming a graph into a series of simplicial complexes, called filtration, and subsequently analyzing the homological properties of these shapes (e.g., the number of connected components, loops, and voids) to gain insights about the graph~\cite{ghrist2008barcodes}. PH has been employed in a variety of graph mining problems, including clustering~\cite{keil2018ego}, node classification~\cite{nguyen2020bot}, graph classification~\cite{rieck2019persistent}, and graph neural networks~\cite{hofer2020graph,horn2021topological,carriere2020perslay}, among others.

When working with graphs, defining persistent homology (PH) while preserving the original graph is a straightforward process, as graphs inherently possess different intrinsic distances and a simplicial complex structure. However, in the case of hypergraphs, preserving topological information without adding or removing structures becomes challenging, as hypergraphs do not necessarily exhibit a simplicial complex structure. The conventional approach for computing the homology of hypergraphs, known as the \textit{simplicial complex closure (SCC)}, involves adding all unobserved sub-edges of a hyperedge. However, this method introduces numerous new hyperedges to the hypergraph, potentially leading to the loss of the original hypergraph structure.



Besides, the work presented in \cite{barycentric} introduces two distinct topological characterizations derived from a hypergraph, which are utilized to investigate hypergraph homology in static hypergraphs through the process of barycentric subdivision applied to the simplicial complex closure. The first characterization, referred to as the \textit{restricted barycentric subdivision (ResBS)}, involves generating the barycentric subdivision of the closure while limiting the subdivision to the hyperedges exclusively. This restriction enables the elimination of hyperedges that are not observed within the hypergraph. On the other hand, the second characterization, known as the \textit{relative barycentric subdivision (RelBS)}, also generates the barycentric subdivision of the closure but additionally employs a quotient operation to remove the subcomplex induced by the missing hyperedges. By employing this quotient, the structure is able to preserve hyperedge interactions while eliminating unobserved hyperedges. Each of these structures emphasizes different higher-order characteristics inherent in hypergraphs while still preserving the original hypergraph structure.

In this paper, we define three distinct filtrations for hypergraph persistent homology using simplicial complex closure, restricted barycentric subdivision, and relative barycentric subdivision of hypergraphs. Our approach involves creating a nested sequence of simplicial complexes for a given hypergraph and defining SCC, ResBS, and RelBS filtrations on each sequence. 

To showcase the effectiveness of hypergraph persistent homology in the field of hypergraph mining, we employ our method to address the hypergraph classification problem. Our approach involves establishing filtrations for a given hypergraph, from which we obtain 0-dimensional and 1-dimensional barcodes for each filtration. By extracting various numerical measures associated with these barcodes, we derive informative hypergraph features.

To assess the performance of our filtrations across different domains within the hypergraph classification problem, we conduct experiments on four real-world hypergraphs: two social networks (Highschool and Primary), one music network (Makam), and one text network (BBC). The results of these experiments demonstrate that our filtration methods exhibit remarkable performance in hypergraph classification, achieving an accuracy rate of approximately 96\%. For comparative analysis, we also utilize graph convolutional networks (GCN) as the baseline method. In this comparison, we project the hypergraphs onto graphs by mapping hyperedges to pairwise edges (i.e., clique expansion) and apply graph convolutional networks (GCN) with global pooling for the purpose of classification. However, our findings indicate that hypergraph persistent homology consistently outperforms the baseline GCN model.

This paper is organized as follows. In Section \ref{sec:rel}, we present the necessary background information. We begin by providing a formal definition of graphs and hypergraphs, followed by an introduction to simplicial complexes and persistent homology. Additionally, we review relevant works on persistent homology in the context of both graphs and hypergraphs. In Section \ref{sec:met}, we delve into the details of the employed hypergraph filtrations and elaborate on our methods for generating hypergraph features based on barcodes. Section \ref{sec:exp} is dedicated to presenting the datasets used in our experiments. We describe how we define the filtration for each dataset and report the obtained classification results on various real-world datasets, comparing them with the results of baseline methods. Finally, we provide concluding remarks in Section \ref{sec:con}.

\section{Preliminaries and Background} \label{sec:rel}

\subsection{Graphs and Hypergraphs}\vspace{-1mm}

\textit{Graphs} are data structures representing relationships between objects \cite{aggarwal2010managing,Cook2006}. They are formed by a set of \textit{vertices} (also called nodes) and a set of \textit{edges} that are connections between pairs of vertices. In a formal definition, a graph $G$ is a pair of sets $G = (V, E)$ where $V$ is the set of vertices and $E \subset V \times V $ is the set of edges of the graph.
There are various types of graphs that represent the differences in the relations between vertices. While in an \textit{undirected graph}, edges link two vertices ${v, w}$ symmetrically, in a \textit{directed graph}, edges link two vertices asymmetrically. If there is a score for the relationship between vertices that could represent the strength of interaction, we can represent this type of relationships or interactions by a \textit{weighted graph}. In a weighted graph, a weight function $W: E \rightarrow \mathbb{R}$ is defined to assign a weight for each edge.

A \textit{hypergraph} $H$ denoted by $H=(V,E=(e_i)_{i \in I})$ on the finite vertex set $V$ is a family $(e_i)_{i \in I}$ ($I$ is a finite set of indexes) of subsets of $V$ called \textit{hyperedges}. Hyperedges can be in different sizes possibly ranging from one vertex $\{v\}\subseteq V$ to the entire vertex set $V$.

\vspace{-2mm}
\subsection{Simplical Complexes} \label{sec:hyperLap}\vspace{-1mm}

A \textit{simplicial complex} is a topological object which is built as a union of points, edges, triangles, tetrahedrons, and higher-dimensional polytopes, i.e. \textit{simplices}.  A 0-simplex is a point, a 1-simplex is two points connected with a line segment, and a 2-simplex is a filled triangle. Here we present the formal definition of a simplicial complex.

\begin{definition}
A simplicial complex $K$ is a finite collection of simplices such that every face of a simplex of $K$ belongs to $K$ and the intersection of any two simplices of $K$ is a common face of both of them.
\end{definition}

To define the topological features of a simplicial complex, we need to understand the concept of a hole in topology. An $n$-dimensional hole can be considered as a void surrounded by $n$-simplices. For example, 0-dimensional holes are connected components, 1-dimensional holes are loops bounded by 1-simplices (i.e., edges), and 2-dimensional holes are voids surrounded by 2-simplices (i.e., triangles). 

Once a graph has a simplicial complex structure, the next step is to identify the presence of any holes within it. This can be accomplished through \textit{simplicial homology}, which allows us to detect and quantify the holes present in a simplicial complex. For a given simplicial complex $X$, simplicial homology associates vector spaces $H_i(X)$ for $i=\{0,1,2,...\}$ where the dimension of $H_i(X)$ gives the number of $i$-dimensional holes. 

\begin{figure}[h!]
    \centering
    \includegraphics[width=.84\textwidth]{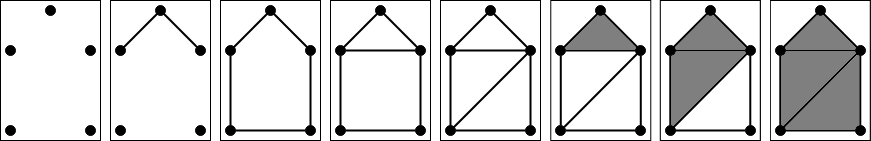}
    \caption{A filtration for $\delta=0,1,2,3,4,5,6,7$ (from left to right), (borrowed from \cite{aktas2019persistence})}
    \label{fig:filter}
\end{figure}

\begin{figure*}[h!]
    \centering
    \begin{subfigure}[t]{0.5\textwidth}
        \centering
        \includegraphics[width=.8\textwidth]{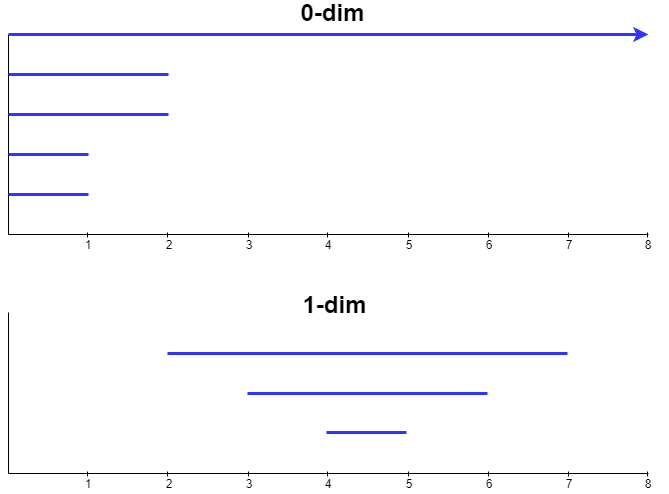}
        \caption{}
    \end{subfigure}%
    ~ 
    \begin{subfigure}[t]{0.5\textwidth}
        \centering
        \includegraphics[width=.7\textwidth]{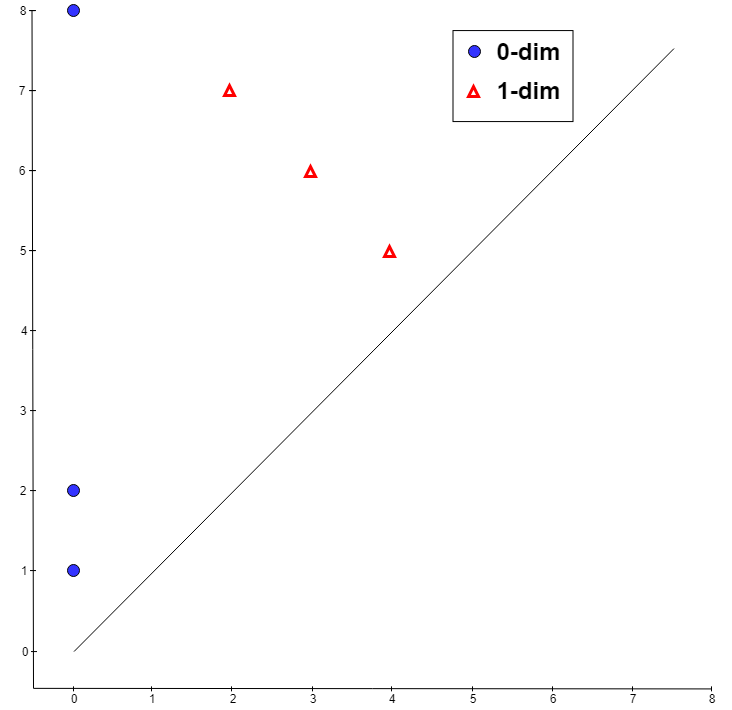}
        \caption{}
    \end{subfigure}
    \caption{Results of the filtration in Figure \ref{fig:filter}, (a) Persistence barcodes for 0- and 1-dim (b) Persistence diagrams for 0- and 1- dim (borrowed from \cite{aktas2019persistence})}
    \label{fig:pd}
\end{figure*}
\vspace{-2mm}

\begin{figure}[h!]
    \centering
    \includegraphics[width=.45\textwidth]{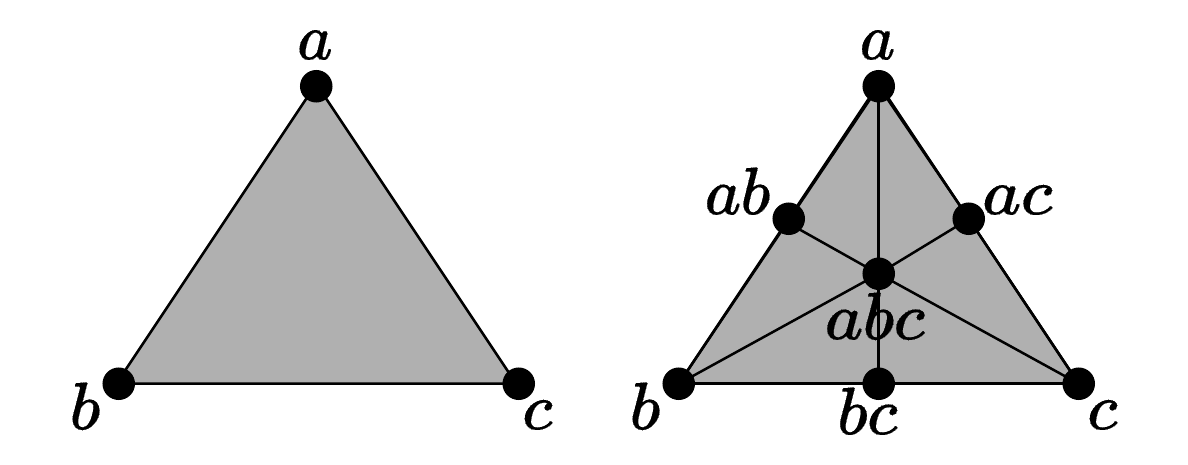}
    \caption{Barycenteric subdivision of the 2-simplex (triangle) on the left.}
    \label{fig:bary}
\end{figure}

On the other hand, the \textit{barycentric subdivision} is an important operation in algebraic topology and computational geometry. It involves dividing each simplex into smaller simplices by adding a new vertex at the center of the simplex and connecting it to each of its subfaces. This process generates a new simplicial complex with more vertices, edges, and higher-dimensional simplices. It allows to refine simplicial complexes while preserving their topological properties. For example, the barycentric subdivision of the triangle in Figure \ref{fig:bary} on the left is given on the right. We replace the 2-simplex with a vertex $abc$ and connect this vertex with its subfaces $a,b,c,ab,ac,bc$. 

\subsection{Persistent homology}\vspace{-1mm}

When analyzing graphs for mining purposes, relying solely on simplicial homology to examine their topological characteristics may not provide sufficient information. For instance, comparing connected graphs based on 0-dimensional homology alone is impractical since both graphs would have the same 0-dimensional homology, which represents the number of connected components. To overcome this limitation, we introduce a family of simplicial complexes derived from a given hypergraph, and observe how the homology evolves within this family. This nested family of simplicial complexes is referred to as a \textit{filtration}, and \textit{persistent homology} is a valuable tool utilized to track the birth and death of these homological features throughout the filtration process. The lifespan of each feature can be represented as an interval, with the initiation and termination points of the interval corresponding to the birth and death of the feature, respectively \cite{edelsbrunner2000topological}.

When examining a dataset using a filtration, \textit{persistence barcodes} capture the intervals that represent the lifespan of homological features. These barcodes offer insights into the duration for which each feature persists, with each bar in the barcode corresponding to a distinct feature. Utilizing persistence barcodes allows us to strike a balance between encoding comprehensive topological information and maintaining computational efficiency.

Alternatively, persistence barcodes can be represented as \textit{persistence diagrams}. These diagrams depict the birth and death times of the features as points (birth, death) in the extended real plane $\bar{\mathbb{R}}^2$. This representation provides a visual understanding of the features' lifespan (see Figure~\ref{fig:filter} and \ref{fig:pd} for an example of a filtration defined on a graph, as well as its corresponding persistent barcodes and diagrams, respectively). In persistence barcodes, longer bars or points further away from the diagonal in persistence diagrams are considered to represent genuine features of the dataset, while others are seen as noise.

\subsection{Related Work} \label{sec:related}

Persistent homology has gained increasing relevance in various graph mining problems in recent years. The applications of persistent homology can be broadly categorized into two main groups: single-graph analysis and multiple-graph analysis. In single-graph analysis, persistent homology is employed to uncover the global structural characteristics of an individual graph. This includes studying the complexity and distribution of strongly connected regions within the graph. For example, researchers have utilized persistent homology to analyze brain graphs, investigating abnormal white matter in maltreated children \cite{chung2013persistent}. In contrast, multiple-graph analysis involves using persistent homology to compare and classify different graphs. For instance, in the work presented in \cite{huang2017persistent}, persistent homology was applied to detect collaboration patterns in collaboration networks. Readers interested in further exploring the applications of persistent homology in graph analysis can find more detailed information in the survey conducted in \cite{aktas2019persistence}.

While there exists a significant body of research on persistent homology (PH) applied to graphs, the exploration of PH in the context of hypergraphs has received comparatively less attention. However, some notable studies have addressed this gap. In \cite{bressan2016embedded} and \cite{ren2020stability}, the authors proposed the construction of embedded homology specifically tailored for hypergraphs. They also demonstrated the stability of PH in the hypergraph setting. Furthermore, Lie et al. utilized embedded homology for hypergraphs to generate molecular descriptors for drug design, while \cite{huang2016persistent} defined a similar filtration approach using clique expansion and leveraged PH to establish a lower bound on higher-order  distances. Additionally, in \cite{tallon2020community}, clique expansion was employed to investigate the topological structures of hypergraphs. More recently, Myers et al. \cite{myers2023topological} analyze the evolution of topological structure in time-evolving hypergraphs. They demonstrate the application of this approach to cyber security and social network datasets, showcasing how the changing topological structure of temporal hypergraphs can provide insights into the underlying dynamics of the systems being studied. In the comprehensive framework for PH methods in hypergraphs presented by Ren in \cite{ren2020persistent}, the existing approaches are summarized and discussed. 

\begin{figure*}[t!]
    \centering
    \begin{subfigure}[t]{0.18\textwidth}
        \centering
        \includegraphics[width=.8\textwidth]{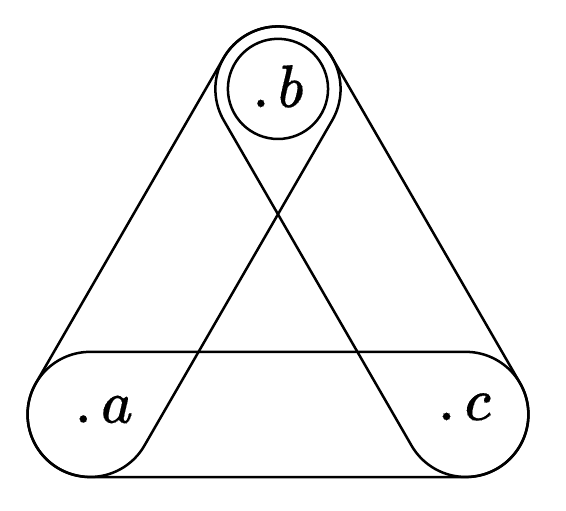}
        \caption{Original hypergraph}
    \end{subfigure}%
    ~ 
    \begin{subfigure}[t]{0.18\textwidth}
        \centering
        \includegraphics[width=.8\textwidth]{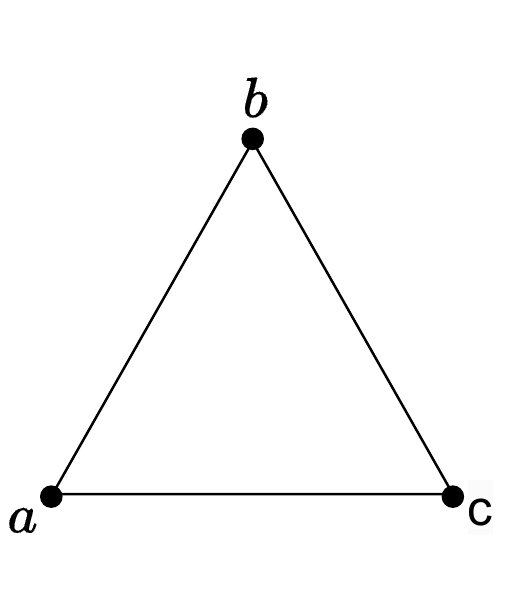}
        \caption{SCC}
    \end{subfigure}
        ~ 
    \begin{subfigure}[t]{0.18\textwidth}
        \centering
        \includegraphics[width=.8\textwidth]{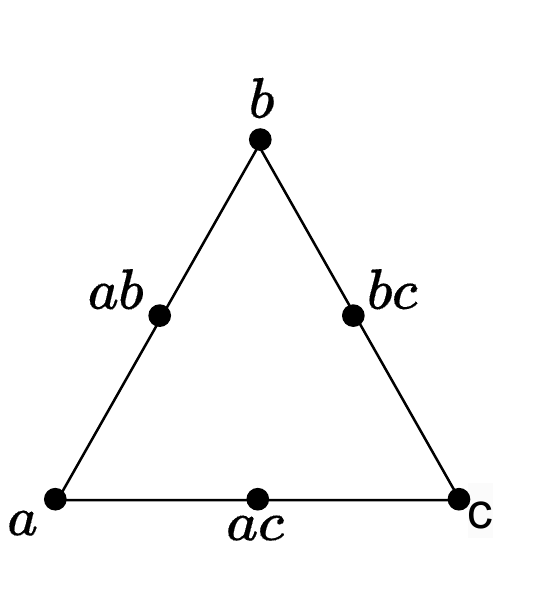}
        \caption{Barycentric subdivision}
    \end{subfigure}
        ~ 
    \begin{subfigure}[t]{0.18\textwidth}
        \centering
        \includegraphics[width=.8\textwidth]{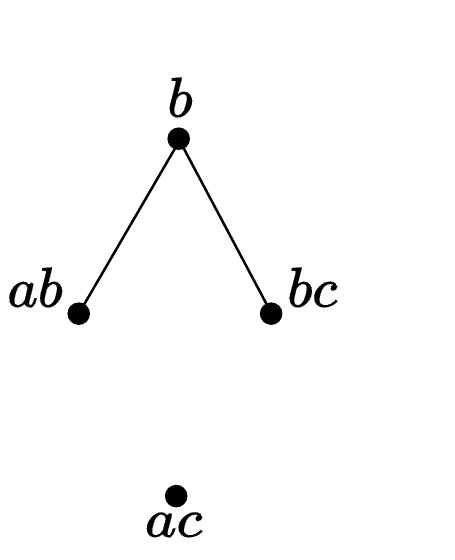}
        \caption{ResBS}
    \end{subfigure}
        ~ 
    \begin{subfigure}[t]{0.18\textwidth}
        \centering
        \includegraphics[width=.8\textwidth]{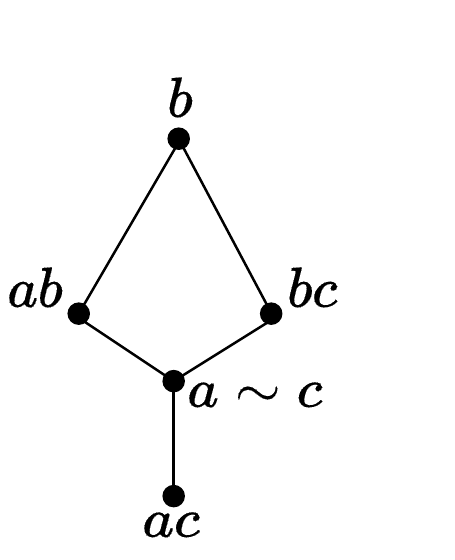}
        \caption{RelBS}
    \end{subfigure}
    \caption{Illustrating various topological characterizations of a hypergraph: (a) The hypergraph with hyperedge set ${b, ab, bc, ac}$. (b) Simplicial complex closure (SCC). (c) Barycentric subdivision. (d) Restricted barycentric subdivision (ResBS). (e) Relative barycentric subdivision (RelBS).}
    \label{fig:bs}
\end{figure*}
\vspace{-2mm}
   
\section{Methods}\label{sec:met}

In this section, we will introduce three hypergraph filtrations that we use for the topological characterization of hypergraphs: simplicial complex closure (SCC), restricted barycentric subdivision (ResBS), and relative barycentric subdivision (RelBS). We further explain how we derive hypergraph features from the persistent diagrams and utilize these features for classification.

\subsection{Hypergraph Filtrations}

\subsubsection{Simplicial complex closure (SCC) filtration} The conventional and widely used filtration for hypergraphs is the \textit{simplicial complex closure} (SCC) filtration. In SCC, we employ the simplicial complex closure of hypergraphs, i.e., we first get the clique expansion of hypergraph and take the clique complex of the expansion. This means that when a hyperedge is added, we include all unobserved sub-edges in the filtration. An illustration of this process is depicted in Figure \ref{fig:bs}.

For the filtration, let $H=\{H_0 \subseteq H_1 \subseteq \cdots H_n\}$ be a nested sequence of hypergraphs. Then, we define the \textit{simplicial complex closure} filtration as follows
$$
\{SCC(H_{i}) \hookrightarrow SCC(H_{j}) \}_{i \leq j}.
$$
Here, $SCC(H_i)$ gives a simplicial complex for all $i \in \{1,...,n\}$, so we get a filtered simplicial complex.

One disadvantage of this filtration is that we add or remove information from the hypergraph, which means we cannot maintain the original hypergraph structure. For instance, removing hyperedge ${b}$ from the hypergraph in Figure \ref{fig:bs}-a still gives the same closure. Nonetheless, this filtration still provides a valuable topological characterization of hypergraphs that can be used in hypergraph analysis.

\subsubsection{Restricted barycentric subdivision filtration} 

For this filtration, we use the restricted barycentric subdivision ($ResBS$) of a hypergraph defined in \cite{barycentric}. For this construction, we first take its simplicial complex closure, then construct the barycentric subdivision of the closure, and finally restrict the subdivision to the hyperedges only. We illustrate how to generate the restricted barycentric subdivision in Figure \ref{fig:bs}.

This subdivision captures different structures on hypergraphs, but intersections between hyperedges are not preserved unless they are present as hyperedges in the hypergraph. For example, in Figure \ref{fig:bs}, the intersection between the hyperedge $ac$ and $bc$ is not preserved in the restricted barycentric subdivision unless $c$ is present as a hyperedge in the hypergraph.  
For the filtration, let $H=\{H_0 \subseteq H_1 \subseteq \cdots H_n\}$ be a nested sequence of hypergraphs. Then, we define the \textit{restricted barycentric subdivision} filtration as follows
$$
\{ResBS(H_{i}) \hookrightarrow ResBS(H_{j}) \}_{i \leq j}.
$$
Here, $ResBS(H_i)$ gives a simplicial complex for all $i \in \{1,...,n\}$ since the restriction of a simplicial complex is also a simplicial complex. Hence, the construction creates a filtered simplicial complex.

\subsubsection{Relative barycentric subdivision filtration} For this filtration, we use the relative barycentric subdivision ($RelBS$) of a hypergraph defined in \cite{barycentric}. For this construction, we again take its simplicial complex closure, then construct the barycentric subdivision of the closure. Next, we quotient out the subcomplex induced by missing hyperedges to preserve the hyperedge intersections. For example, in Figure \ref{fig:bs}, when we quotient out the missing hyperedges ${a}$ and ${c}$, we can preserve the intersection between $ac$ and $bc$ through the vertex $c$ (or, similarly for the intersection between $ac$ and $ab$ through the vertex $a$). This filtration creates one-dimensional loops when we quotient out hyperedges.

One potential issue arises from the fact that the quotient of a simplicial complex may not necessarily be a simplicial complex itself. Therefore, it is essential to ensure that we have a filtered simplicial complex before computing persistent homology. In our approach, constructing the barycentric subdivision of the simplicial complex closure represents each simplex in the closure as a vertex in the subdivision. When quotienting out the subcomplex induced by missing hyperedges, we only identify and glue the corresponding vertices with their edges in the subdivision, discarding any higher-dimensional simplices. There is a specific case where this construction does not result in a simplicial complex, which occurs when two edges are created between two vertices. To address this issue, we remove one of the edges between these vertices. This modification ensures the creation of a simplicial complex while preserving the desired homology properties. 

For the filtration, let $H=\{H_0 \subseteq H_1 \subseteq \cdots H_n\}$ be a nested sequence of hypergraphs. Then, we define the \textit{restricted barycentric subdivision} filtration as follows
$$
\{RelBS(H_{i}) \hookrightarrow RelBS(H_{j}) \}_{i \leq j}.
$$
Incorporating the modification described in the previous paragraph within the relative barycentric subdivision, we ensure that $RelBS(H_i)$ forms a simplicial complex for all $i \in {1,...,n}$.

\subsection{Feature Extraction} \label{sec:feat}
Persistent barcodes contain valuable topological information about hypergraphs. However, since barcodes consist of intervals rather than numerical features, they cannot be directly utilized in classification models. Therefore, we define five numerical features for each dimension based on the barcodes. The first feature represents the number of bars in each barcode, which corresponds to the count of connected components for 0-dimensional barcodes and cycles for 1-dimensional barcodes. This count plays a crucial role in distinguishing between different hypergraph classes. Additionally, we introduce four more features inspired by \cite{adcock2013ring} as follows
\[
\sum_i x_i(y_i-x_i),
\]
\[
\sum_i (y_{max} - y_i)(y_i - x_i),
\]
\[
\sum_i x_i^2(y_i-x_i)^4,
\]
\[
\sum_i (y_{max} - y_i)^2(y_i-x_i)^4.
\]
In these formulas, $x_i$ and $y_i$ represent the start and endpoint of each interval $i$ in a barcode, and $y_{max}$ denotes the maximum $y_i$ value within the barcode. The first two features incorporate information from all bars, including their lengths and endpoints, while the latter two features focus on the arrangement of longer bars. To ensure proper definition of these features, we assign the filtration value of infinite bars to their maximum value, represented as $y_{max}$.

\subsection{Classification}
After obtaining the feature vector for each hypergraph, we utilize the Random Forest classification algorithm to build our prediction model. Random Forest is a well-known algorithm that constructs multiple decision trees and combines them to achieve a more accurate and robust prediction. To evaluate our method, we employ the 10-fold cross-validation process, where the dataset is randomly divided into 10 approximately equal-sized folds. Each fold is treated as a validation set, while the remaining folds are used for training the model.

\section{Experiments}\label{sec:exp}
In this section, we first provide an overview of the experimental setup, including information about the real-world datasets used and the baseline models employed for comparison. Secondly, we detail the process of generating filtered simplicial complexes for each dataset, along with the visualization of the resulting persistent diagrams for a representative hypergraph from the datasets. Lastly, we present the experimental results, which involve comparing the performance of the proposed filtration models among themselves and also against the baseline models. 

\begin{table}[h!]
\centering

\caption{Properties of the graph datasets used in experiments. $|H|$ is the number of hypergraphs, $|E|$ is the total number of hyperedges, and $k_{max}$ is the maximum hyperedge size.} \vspace{3mm}
\begin{tabular}{|c|c|c|c|}

\hline \cellcolor{gray!60}\textbf{Datasets}  & \cellcolor{gray!60} $\mathbf{|H|}$  &\cellcolor{gray!60} $\mathbf{|E|}$ & \cellcolor{gray!60}$\mathbf{k}_{max}$ \\
 \hline \hline  \cellcolor{gray!25}\textbf{Highschool} & 327 & 7,818 & 5 \\	
\hline \cellcolor{gray!25}\textbf{Primary} & 242 & 12,704 & 5\\	
\hline \cellcolor{gray!25}\textbf{Makam} & 832 & 245,534 & 87 \\
\hline
\cellcolor{gray!25}\textbf{BBC} & 2225 & 41,403 & 159\\
\hline 

\end{tabular}
\label{table:data}
\end{table}

\subsection{Setup}
\subsubsection{Datasets}
In our experimental setup, we assess the effectiveness of the proposed filtrations in hypergraph classification using four undirected real-world networks. The network statistics for each dataset are provided in Table~\ref{table:data}, where $|H|$ represents the number of hypergraphs, $|E|$ denotes the total number of hyperedges, and $k_{max}$ indicates the maximum hyperedge size. The datasets used in our evaluation include: (1) Highschool~\cite{Stehl-2011-contact,chodrow2021hypergraph}, which is a network of high school students in Marseilles, France. Each vertex represents a student, and a hyperedge represents a set of students in close contact with each other. We generate ego networks for each student to obtain the hypergraph set. (2) Primary school~\cite{Mastrandrea-2015-contact,chodrow2021hypergraph}, which is a network of primary school students and teachers. Each vertex corresponds to a student or a teacher, and a hyperedge denotes a set of students and/or teachers in close contact. Similarly, we generate ego networks for each student to obtain the hypergraph set. (3) Makam~\cite{aktas2019classification, karaosmanouglu2012turkish}, which consists of Turkish Makam music songs. For each song, we create a hypergraph where the vertices represent notes, and each syllable forms a hyperedge comprising the corresponding notes. (4) BBC~\cite{greene2006practical}, which is a dataset comprising BBC news articles. For each article, we construct a hypergraph with words as vertices, and each sentence forms a hyperedge. The datasets can be accessed at \url{https://sites.google.com/view/mehmetaktas/datasets}.

\begin{figure}[t!]
 
\centering
\includegraphics[width=\textwidth]{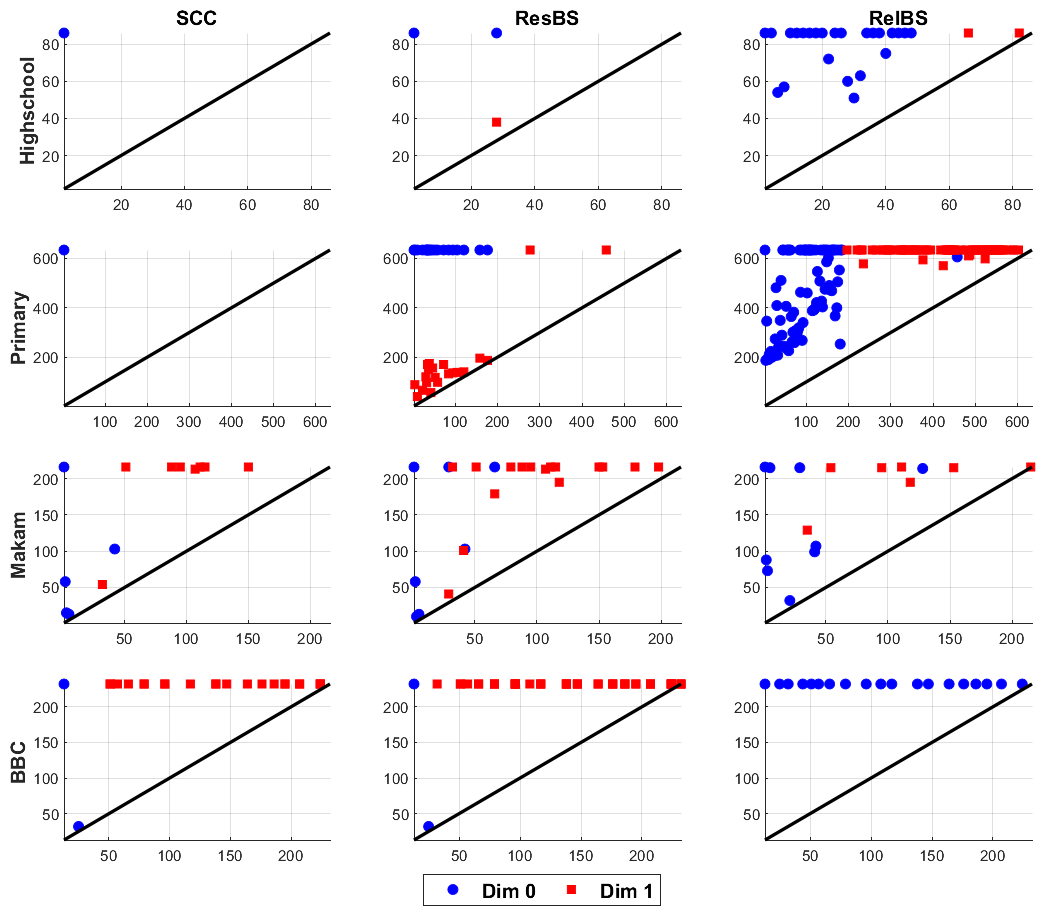}
\caption{0- and 1-dimensional persistent diagrams of one hypergraph within the datasets in Table \ref{table:data} using SCC, ResBS, and RelBS filtrations.}
\label{fig:filt}
   \end{figure}
   
\subsubsection{Filtrations}
To establish the filtration for the given datasets, we follow specific procedures based on the characteristics of each dataset. For the Highschool and Primary datasets, which capture interactions recorded by wearable sensors in schools at a 20-second resolution, we incorporate the recorded time information into the filtration process. In the case of the Makam music dataset, hyperedges are added according to their order of appearance in the musical note sheet for each song. Similarly, for the BBC text network dataset, hyperedges are added based on their sequential occurrence in the corresponding news articles. As a computational efficiency measure, we only include hyperedges with a size less than 25 in all datasets during the filtration process for our experiments. Consequently, we obtain a nested sequence of simplicial complexes for each dataset, which forms the basis for defining persistent homology.

In Figure \ref{fig:filt}, we display 0- and 1-dimensional persistent diagrams of one hypergraph from each dataset in Table \ref{table:data} using SCC, ResBS, and RelBS filtrations. The figure indicates that each hypergraph has unique topological characteristics. Highschool and Primary are networks of school interactions. Due to the rapid interaction between students in schools, the corresponding SCC persistent diagrams have only one bar for 0-dimension and no bar for 1-dimension. Notably, these networks lack 1-dimensional cycles. Conversely, the Primary network has numerous bars for both dimensions, particularly with ResBS and RelBS filtrations. This could be attributed to the inclusion of teachers in the network. 

In contrast, Makam hypergraphs, which represent musical chords, have multiple bars for all filtrations in both dimensions. This is because of the non-linear complex structure of musical chords. Lastly, the BBC hypergraph, where each sentence forms a hyperedge, has persistent diagrams with bars extending to infinity due to disconnected components resulting from different words in hyperedges. It is noteworthy that there is no bar for 1-dimensional in the ResBS filtration of the BBC dataset.

\subsubsection{Baseline Models} \label{sec:base}

We compare our hypergraph filtrations with Graph Neural Networks (GNNs) as the baseline model. GNNs have recently demonstrated state-of-the-art performance for graph-related tasks such as node classification, graph classification, and community detection. To compare the GNN model with our model for the graph classification task, we use Graph Convolutional Network (GCN) \cite{kipf2017semisupervised}, which is one of the most popular GNN models. We follow the global pooling architecture \cite{zhang2018end} to learn graph-level representation. In this architecture, there are three message-passing layers where we use the GCN model for message-passing. Each message-passing layer updates the node representation. The output of each message-passing layer is concatenated and given to a global pooling layer \cite{zhang2018end} to obtain the graph-level representation. A readout layer \cite{cangea2018towards} is used to process the output of the pooling layer, and an MLP layer is used for graph classification.

\subsection{Results}
This section presents the results of our experiments on each dataset. We begin by comparing the performance of the proposed filtration models. Then, we compare the performance of these models with the GNN model described in the previous section. For both experiments, we use 0-dimensional and 1-dimensional persistent barcode features. 

\begin{figure}[h!]
 
\centering
\includegraphics[width=.6\textwidth]{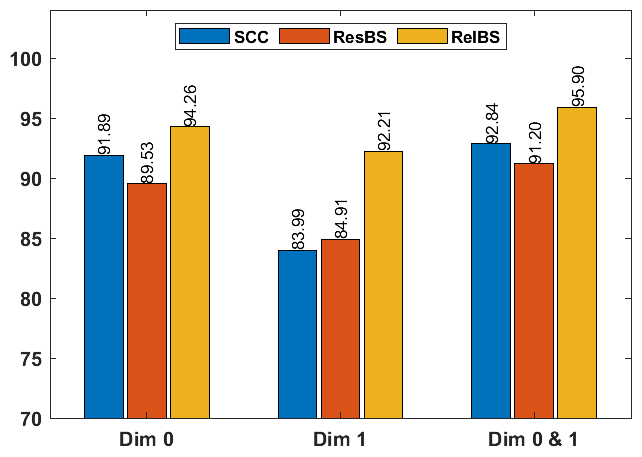}
\caption{Performance comparison (accuracy) between our filtrations, ResBS and RelBS, with SCC.}
\label{fig:acc}
   \end{figure}

We first create the 0-dimensional and 1-dimensional barcodes for each hypergraph using SCC, ResBS, and RelBS filtrations. Then, we obtain the features for each barcode as it is outlined in Section~\ref{sec:feat}. In our experiment, besides using the features of 0-dimensional and 1-dimensional barcodes, we also concatenate these features for each barcode to utilize them together. The results are available in Figure \ref{fig:acc} and Table \ref{table:compsin}. 

Figure \ref{fig:acc} illustrates the classification accuracy results, which indicate that the RelBS filtration method outperforms other filtration methods for all dimensions. In terms of SCC and ResBS filtrations, SCC achieves better performance for 0-dimensional features, whereas ResBS filtration performs better for 1-dimensional features. The reason behind this is that the Primary dataset's 1-dimensional barcodes are empty for the SCC filtration since there is no 1-dimensional cycle in the filtration.

Moreover, the accuracy results improve when 0-dimensional and 1-dimensional features are used together. This is expected as different dimensions capture different topological characteristics of hypergraphs. Hence, utilizing these features together in the classifier yields better results.

\begin{table*}

\centering{
\caption{Performance comparison between the proposed filtration models based on precision and recall measures. The best result for each dataset in each measure is typed in bold.}
\label{table:compsin}
\normalsize
\begin{tabular}{|c | c c c | c c c | }
 \hline

    \multirow{2}{*}{\textbf{Dim 0}} & \multicolumn{3}{|c|}{\cellcolor{gray!60}\textbf{Precision}} & \multicolumn{3}{|c|}{\cellcolor{gray!60}\textbf{Recall}} \\ \cline{2-7}
 & \cellcolor{gray!25}\textbf{SCC} & \cellcolor{gray!25}\textbf{ResBS} & \cellcolor{gray!25}\textbf{RelBS} & \cellcolor{gray!25}\textbf{SCC} & \cellcolor{gray!25}\textbf{ResBS} & \cellcolor{gray!25}\textbf{RelBS} \\ \hline
 
 \cellcolor{gray!25}\textbf{Highschool} & 76.30 & 76.10 & \textbf{79.30} & 78.60 & 76.10 & \textbf{85.30} \\ \hline
 
  \cellcolor{gray!25} \textbf{Primary} & 69.30 & 64.80 & \textbf{78.10} & 69.00 & 65.30 & \textbf{70.70} \\ \hline
   \cellcolor{gray!25} \textbf{Makam} & 94.90 & 91.30 & \textbf{96.10} & 88.20 & 84.10 & \textbf{94.20} \\ \hline
   \cellcolor{gray!25} \textbf{BBC} & 95.80 & 93.80 & \textbf{97.80} & 98.20 & 96.70 & \textbf{98.50} \\ \hline 
\hline
    \multirow{2}{*}{\textbf{Dim 1}} & \multicolumn{3}{|c|}{\cellcolor{gray!60}\textbf{Precision}} & \multicolumn{3}{|c|}{\cellcolor{gray!60}\textbf{Recall}} \\ \cline{2-7}
 & \cellcolor{gray!25}\textbf{SCC} & \cellcolor{gray!25}\textbf{ResBS} & \cellcolor{gray!25}\textbf{RelBS} & \cellcolor{gray!25}\textbf{SCC} & \cellcolor{gray!25}\textbf{ResBS} & \cellcolor{gray!25}\textbf{RelBS} \\ \hline
 
 \cellcolor{gray!25}\textbf{Highschool} & 46.60 & 62.50 & \textbf{73.30} & \textbf{100} & 92.70 & 56.30 \\ \hline
 
  \cellcolor{gray!25} \textbf{Primary} & NA & \textbf{79.20} & 78.50 & 0.00 & 39.30 & \textbf{66.50} \\ \hline
   \cellcolor{gray!25} \textbf{Makam} & 90.70 & 82.20 & \textbf{98.60} & 77.70 & 74.50 & \textbf{94.60} \\ \hline
   \cellcolor{gray!25} \textbf{BBC} & \textbf{94.40} & 91.40 & 93.30 & 93.80 & 93.20 & \textbf{100} \\ \hline 
\hline
    \multirow{2}{*}{\textbf{Dim 0 \& 1}} & \multicolumn{3}{|c|}{\cellcolor{gray!60}\textbf{Precision}} & \multicolumn{3}{|c|}{\cellcolor{gray!60}\textbf{Recall}} \\ \cline{2-7}
 & \cellcolor{gray!25}\textbf{SCC} & \cellcolor{gray!25}\textbf{ResBS} & \cellcolor{gray!25}\textbf{RelBS} & \cellcolor{gray!25}\textbf{SCC} & \cellcolor{gray!25}\textbf{ResBS} & \cellcolor{gray!25}\textbf{RelBS} \\ \hline
 
 \cellcolor{gray!25}\textbf{Highschool} & 78.30 & 78.50 & \textbf{79.80} & 79.50 & 79.20 & \textbf{86.90} \\ \hline
 
  \cellcolor{gray!25} \textbf{Primary} & 71.40 & 71.50 & \textbf{81.70} & 70.20 & \textbf{70.70} & 70.20 \\ \hline
   \cellcolor{gray!25} \textbf{Makam} & 95.70 & 93.70 & \textbf{99.30} & 91.00 & 86.00 & \textbf{97.10} \\ \hline
   \cellcolor{gray!25} \textbf{BBC} & 96.50 & 94.50 & \textbf{98.80} & 98.40 & 97.60 & \textbf{99.90} \\ \hline 

    \end{tabular}}
 \vspace{-0.25cm}
\end{table*}

In Table \ref{table:compsin}, we present the precision and recall results for each dataset in addition to accuracy. For Dim 0, RelBS filtration outperforms SCC and ResBS filtration in terms of both precision and recall. However, for Dim 1, there is no single filtration that outperforms the others consistently.

In the Highschool dataset, RelBS filtration has the best precision with 73.30\%, while SCC has a perfect recall score. In the Primary dataset, ResBS filtration has the best precision with 79.20\%, while RelBS filtration has the best recall with 66.50\%. In the Makam dataset, ResBS filtration has the best precision and recall. In the BBC dataset, SCC filtration has the best precision with 94.40\%, and ResBS filtration has a perfect recall score. Finally, for Dim 0 \& 1, RelBS filtration has the best precision and recall except for the recall score in the Primary dataset.
   
\begin{figure}[h!]
 
\centering
\includegraphics[width=.4\textwidth]{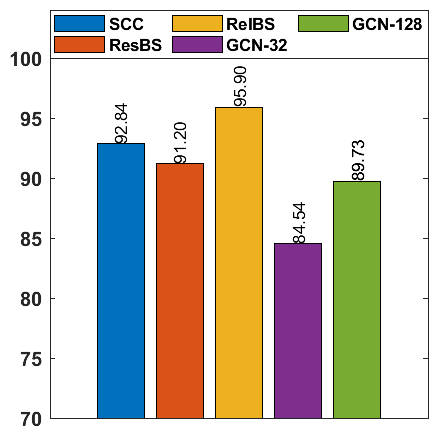}
\caption{Performance comparison (accuracy) between persistent homology methods and graph convolutional network (GCN). We use RelBS filtration as the persistent homology method since it provides the best accuracy. We further select the hidden size parameter in GCN as 32 and 128 to see the effect of the feature size.}
\label{fig:comp}
   \end{figure}

In addition to our proposed filtrations, we also compare them with graph neural networks (GNNs). To enable this comparison, we first project each hypergraph into a graph representation using clique expansion, which allows us to treat hypergraph classification as graph classification. We then utilize a graph convolutional network (GCN) model with global pooling as the baseline for our evaluation. To examine the influence of feature sizes, we experiment with hidden sizes of 32 and 128 in the GCN model. To ensure a fair and rigorous comparison, we employ a 10-fold cross-validation process to assess the accuracy of the classification. The results of this comparison can be seen in Figure \ref{fig:comp}.

The figure shows that all filtration models outperform the GCN models selected as the baseline. This is because the GCN models are applied on projected graphs, which means that they overlook higher-order structures present within each hypergraph, leading to an information loss. On the other hand, the filtration models we propose are capable of capturing and leveraging higher-order topological structures within hypergraphs, which results in improved performance compared to the baseline models.

\section{Conclusion}\label{sec:con}
This study presents three distinct filtrations for hypergraph persistent homology, emphasizing various higher-order network structures. The effectiveness of these methods is evaluated through hypergraph classification tasks on four real-world networks from social, music, and text domains. Comparative analysis against a graph neural network (GNN) approach applied to projected graphs demonstrates that the proposed persistent homology filtrations outperform GNN. These findings underscore the potential of utilizing persistent homology in addressing hypergraph mining challenges.

\section*{Acknowledgments}
We would like to thank Emilie Purvine for her feedback about the ideas in the article and also for her wonderful talk ``Homology of Graphs and Hypergraphs" at  the workshop on ``Topological Data Analysis - Theory and Applications" supported by the Tutte Institute and Western University in May 1-3, 2021~\cite{barycentric}\footnote{\url{https://www.youtube.com/watch?v=XeNBysFcwOw}}, that motivated us to initiate this project.

\bibliographystyle{unsrt}  
\bibliography{references}

\end{document}